\documentclass[11pt]{amsart}

\usepackage{amsmath,amssymb,latexsym,amsthm,enumerate}
\usepackage{verbatim}
\usepackage{graphicx} 
\usepackage{hyperref} 

\usepackage{xcolor}

\def\NN{\mathbb{N}}

\def\RR{\mathbb{R}}

\def\PP{\mathbb{P}}

\def\1{\mathbf{1}}
\def\0{\mathbf{0}}

\def\I{\mathbf{I}}

\newcommand{\te}{\widehat{e}}

\newtheorem{Thm}{Theorem}
\newtheorem{Cor}[Thm]{Corollary}
\newtheorem{Lemma}[Thm]{Lemma}
\newtheorem{Prop}[Thm]{Proposition}
\theoremstyle{definition}

\newtheorem{As}{Assumption}

\def\BEN{\begin{enumerate}}  \def\BI{\begin{itemize}}
\def\EEN{\end{enumerate}}   \def\EI{\end{itemize}}

\def\mbb{\mathbb} \def\mbf{\mathbf} \def\mrm{\mathrm}
\def\mc{\mathcal}  \def\ovl{\overline}

\def\le{\left}
\def\ri{\right}
\def\te#1{\mathrm{e}^{#1}}   

\def\WH{\widehat}

\def\g{\gamma}  \def\d{\delta}   
\def\ve{\varepsilon}
\def\e{\epsilon}
 
     \def\s{\sigma}

  \def\td{\text{\rm d}}
 
\numberwithin{equation}{section}

\newcommand{\exit}{{\mbox{\, \vspace{3mm}}}
\hfill\mbox{$\square$}}

\begin{document}

\title[Path functionals of the reflected process]{Joint asymptotic 
distribution of certain path functionals of the reflected process}
\thanks{{\it Acknowledgement:} 
MP's research was supported in part by the NWO-STAR cluster.}

\author{Aleksandar Mijatovi\'{c}}
\address{Department of Mathematics, Imperial College London}
\email{a.mijatovic@imperial.ac.uk}
\author{Martijn Pistorius}
\address{Department of Mathematics, Imperial College London}
\email{m.pistorius@imperial.ac.uk}



\begin{abstract}
Let 
$\tau(x)$ 
be the first time the reflected process
$Y$
of a L\'evy processes
$X$
crosses 
$x>0$. 
The main aim of the paper is to investigate the asymptotic 
dependence of the path functionals:
$Y(t) = X(t) - \inf_{0\leq s\leq t}X(s)$,
$M(t,x)=\sup_{0\leq s\leq t}Y(s)-x$
and
$Z(x)=Y(\tau(x))-x$.
We prove that 
under Cram\'{e}r's condition 
on 
$X(1)$,
the functionals $Y(t)$,  $M(t,y)$ and $Z(x+y)$ are asymptotically
independent as 
$\min\{t,y,x\}\to\infty$.  
We also characterise
the law of the limiting overshoot $Z(\infty)$ 
of the reflected process.
If,
as
$\min\{t,x\}\to\infty$,  
the quantity 
$t\te{-\gamma x}$
has a positive limit
($\gamma$ denotes the Cram\'er coefficient),
our results together with the theorem of 
Doney \& Maller~\cite{DoneyMaller}
imply the existence and the explicit form
of the joint weak limit 
$(Y(\infty), M(\infty),Z(\infty))$.
\end{abstract}

\maketitle

\section{Introduction}
\label{sec:Intro}

The reflected process 
$Y$
of a L\'{e}vy process 
$X$
is a strong Markov process on $\RR_+\doteq[0,\infty)$
equal to 
$X$
reflected at its running
infimum. The reflected process 
is of great importance in many areas
of probability, ranging from the fluctuation
theory for L\'evy processes (e.g.~\cite[Ch.~VI]{Bertoin} and the
references therein)
to  mathematical statistics 
(e.g.~\cite{Mous,Shiryaev}, CUSUM method of cumulative sum), 
queueing theory (e.g.~\cite{AsmussenAPQ, Prabhu}), mathematical
finance (e.g.~\cite{had, MiPi}, drawdown as risk measure), 
mathematical genetics (e.g.~\cite{KarlinDembo} and references therein)
and many more. 
The aim of this paper is to study the asymptotic dependence 
and weak limiting behaviour 
of the functionals of the reflected process $Y$:
\begin{equation}
\label{eq:Levy_Overshoot}
Y(t)\doteq\, X(t) - \inf_{0\leq s\leq t}X(s),\qquad 
M({t,x}) \doteq Y^*(t) - x,\qquad
Z(x)\doteq Y({\tau(x)}) - x,
\end{equation}
where
$t, x \in \mbb \RR_+$.
Here 
$\tau(x)$
and 
$Y^*(t)$
denote the first entry time 
of 
$Y$
into the interval
$(x,\infty)$
and the supremum up to time $t$
of 
the reflected process
respectively,  
\begin{equation*}
\tau(x) \,\doteq\, \inf\{t\ge 0: Y(t) > x\} \quad\text{ ($\inf\emptyset\doteq\infty$)},\qquad
Y^*(t)\doteq\sup_{0\leq s\leq t}Y(s).
\end{equation*}
The main result of the paper identifies a condition
on the L\'evy measure of $X$, under which the triplet
in~\eqref{eq:Levy_Overshoot}
is essentially asymptotically independent in the following sense.
A family of random vectors 
$\{(U_z^1, \ldots, U^d_z)\}_{z\in\RR_+^l}$
on a given probability space,
where
$d,l\in\NN$,
is \textit{asymptotically independent} if
the joint CDF is asymptotically equal to a product 
of the CDFs of the components: i.e. 
for any
$a_i\in(-\infty,\infty]$, $i=1,\ldots,d$,  it holds
(denote 
$a\wedge b\doteq\min\{a,b\}$)
\begin{equation*}
P(U_z^1 \leq a_1, \ldots, U_z^d\leq a_d) = \prod_{i=1}^d P(U^i_z\leq a_i)
 + o(1) 
\qquad
\text{as $z_1\wedge\ldots\wedge z_l\to\infty$.}
\end{equation*}

\noindent Under Assumption~\ref{A2}, which is assumed throughout the paper, our main result, Theorem~\ref{thm},
holds.

\begin{As}\label{A2} 
The mean of $X(1)$ is finite,
\textit{Cram\'{e}r's condition},
$E[\te{\gamma X(1)}]=1$
for $\gamma>0$,
holds,
$E[\te{\gamma X(1)}|X(1)|]<\infty$
and
either 
the 
L\'evy measure of
$X$
is non-lattice or
$0$
is regular for
$(0,\infty)$.
\end{As}

\begin{Thm}\label{thm}
The triplet
$\{(Y(t), Z(x+y), M({t,x}))\}_{t,x,y\in\mbb R_+}$
is asymptotically independent
and the weak limit 
$Z(x) \stackrel{\mc D}{\longrightarrow} Z(\infty)$,
as $x\to\infty$,
holds,
where 
$\phi$
is the Laplace exponent of the increasing ladder-height
process of $X$:\footnote{Note that the Cram\'{e}r condition implies $E[X(1)]<0$ 
and hence
$\phi(0)>0$
(see~\eqref{eq:phi}
for definition of
$\phi$
and Section~\ref{subsec:Setting}
for more details on ladder processes), 
making the formula in~\eqref{eq:zinf} well defined.}
\begin{equation}\label{eq:zinf}
E[\te{-v Z({\infty})}] = \frac{\gamma}{\gamma + v} \cdot \frac{\phi(v)}{\phi(0)}
\qquad\text{for all}\quad v\in\RR_+.
\end{equation}
\end{Thm}

\begin{figure}[t]
\label{fig:3_functionals}
\input{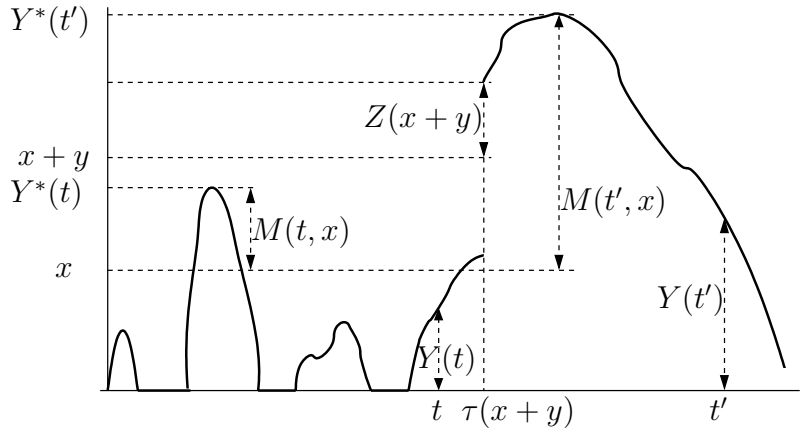}
\caption{\footnotesize{This schematic figure of a
path of 
$Y$
depicts the values of the three functionals in~\eqref{eq:Levy_Overshoot}
at times
$t$
and
$t'$,
before and after the reflected process crosses the
level
$x+y$.
It is intuitively clear that, in general, 
$M(t,x)$,
$Z(x+y)$
and
$Y(t)$
cannot be independent for fixed
$t,x,y>0$.}}
\end{figure}

The asymptotic independence in Theorem~\ref{thm} should be contrasted to 
the intuition that the functionals 
$Z(x+y)$
and
$M(t,x)$
(and hence the triplet in Theorem~\ref{thm})
are unlikely to be asymptotically independent
if the L\'evy measure of 
$X$
is heavy-tailed (e.g. if its tail function is 
regularly varying at infinity).
Intuitively, in this case  
the asymptotic behaviour of the functionals 
is governed by infrequent but very large jumps 
that determine the values of 
$Z(x+y)$
and
$M(t,x)$
simultaneously.
This is analogous to the behaviour of the path
at time
$t'$
in Figure~\ref{fig:3_functionals}.
In fact, 
in contrast to the heavy-tailed case, 
the intuitive reason for the 
asymptotic independence under As.~\ref{A2}
is closely related to the following assertion: 
the likelihood of a single excursion of 
$Y$
straddling both the running time
$t$
and the first-passage time
$\tau(x+y)$,
as depicted in Figure~\ref{fig:3_functionals},
tends to zero
(see remarks following Lemma~\ref{lem:straddle}, Section~\ref{sec:proof2}, 
for a more detailed intuitive explanation of this phenomenon
under As.~\ref{A2}).
This fact will be used to establish the asymptotic independence of 
$M(t,x)$,
$Z(x+y)$
and
$Y(t)$.

It is not hard to see that in general
the functionals
$M(t,x)$
and
$Z(x)$
are not asymptotically independent as
$t\wedge x\to\infty$.
We show (see the remark following 
Lemma~\ref{lem:ind}) that the probability 
that both random variables occur during the 
same excursion of the reflected process
$Y$ away from zero,
as is the case at time
$t'$
in Figure~\ref{fig:3_functionals},
may not decay to zero
under As.~\ref{A2}
However, 
by Theorem~\ref{thm},
the pair 
$\{(M(t,x),Z(x+y))\}_{t,x,y\in\mbb R_+}$ 
is asymptotically independent as 
$t\wedge x\wedge y\to\infty$.
Hence,
for any 
$\alpha>1$
so are the variables
$M(t,x)$ and $Z(\alpha x)$
as
$t\wedge x\to\infty$,



While it may appear intuitively clear that the overshoot 
$Z(x)$
of a \textit{high} level 
$x$
does not occur frequently \textit{during} the excursion straddling
$t'$
but \textit{before} time 
$t'$
(as depicted in Figure~\ref{fig:3_functionals})---indeed, this excursion
still has time to run and reach higher levels, while previous excursions, 
which have concluded their runs, are more likely to have got to the level 
$x$---it does not seems immediately obvious how to make such heuristic arguments precise, 
particularly since in Theorem~\ref{thm} 
the level 
$x$
is allowed to go to infinity arbitrarily slowly compared to the running time
$t$.
One of the contributions of the paper is to establish rigorously the 
asymptotic independence of
$Z(x)$
and
$Y(t)$
(cf. Lemma~\ref{lem:straddle} and remarks that follow).

The definition of asymptotic independence of the 
functionals in~\eqref{eq:Levy_Overshoot}
requires an approximate factorisation of the joint CDF without specifying
the rate of divergence or the mutual dependence of
$t,x,y$
as they tend to infinity. 
Hence the asymptotic independence 
in Theorem~\ref{thm} does not require the existence 
of the weak limit of the functionals.
It is clear however that the most interesting 
application of Theorem~\ref{thm}
is precisely in the case when,
for each of the functionals, 
such a limit exist. 
A result of Doney \& Maller~\cite[Thm~1]{DoneyMaller}
implies that 
$M(t,x)$
converges weakly to a Gumbel distribution 
if the quantity 
$t\te{-\gamma x}$
tends to a positive constant 
as
$t\uparrow\infty$
(see~\cite[Ch.~3]{ExtremeEvents_Book} for the Gumbel 
distribution and Appendix~\ref{subsec:Cor}, equation~\eqref{eq:Ystar},
for a simple derivation of the limit law 
$M(\infty)$
from~\cite[Thm~1]{DoneyMaller}).
Cram\'er's condition implies that 
$X$
tends to 
$-\infty$
almost surely,
and hence, by the classical time reversal argument,
the reflected process
$Y$
has a stationary distribution
$Y(\infty)$
equal to the law of the ultimate supremum 
$\sup_{t\geq0}X(t)$.
The following corollary of Theorem~\ref{thm}
describes explicitly the various limit laws. 

\begin{Cor}
\label{cor:4}
(i) The weak limit of the random vector 
$(Y(t),Z(x))$,
as
$x\wedge t\to\infty$, exists and the law 
$(Y(\infty),Z(\infty))$
is determined by the joint Laplace transform 
\begin{equation*}
E[\exp\left(-u Y(\infty) -v Z(\infty)\right)] = \frac{\gamma}{\gamma + v} \cdot
\frac{\phi(v)}{\phi(u)}
\qquad\text{for all}\quad u,v\in\RR_+.
\end{equation*}

\noindent (ii) Let
$m\doteq\lim_{u\to\infty}\phi(u)/u$
and
$\nu_H$
be the L\'evy measure 
of the Laplace exponent
$\phi$
with the tail function
$\ovl \nu_{H}(x)\doteq\nu_H((x,\infty))$, 
$x>0$.
Then the law of the asymptotic overshoot 
$Z(\infty)$
is given by:
\begin{equation*}
P(Z(\infty)>x) =  \frac{\gamma}{\phi(0)}
\te{-\gamma x}
\int_x^\infty \te{\gamma y}\, \ovl\nu_{H}(y)\,\td y,
\quad\text{$x\in[0,\infty)$,}
\quad\text{and}\quad
P(Z(\infty)=0) = \frac{\gamma}{\phi(0)}m.
\end{equation*}
In particular,
$Z(\infty)$
is a continuous random variable except possibly at 
the origin.

\noindent (iii) Assume that 
$\lim_{t\uparrow\infty}t\te{-\gamma x}= \lambda$ 
for some $\lambda>0$.
Then 
$(Y(t),Z(x+y),M(t,x))$
converges weakly 
as $t\wedge y\to\infty$
and the joint limit law 
$(Y(\infty), Z(\infty), M(\infty))$
is given by the Fourier-Laplace transform: 
\begin{equation*}
E\left[\exp\left(-u Y(\infty) -v Z(\infty) + \mrm i\beta M(\infty)\right)\right] =
\frac{\gamma}{\gamma + v} \cdot \frac{\phi(v)}{\phi(u)} \cdot
\Gamma\le(1 - \frac{\mrm i\beta}{\gamma}\ri)\cdot\exp\left[\mrm i \beta
\gamma^{-1} \log \left(\lambda \ell C_\gamma \WH\phi(\gamma)\right)\right]
\end{equation*}
for all
$u,v\in\RR_+$,
$\beta\in\RR$,
where 
$\WH\phi$ 
is the Laplace exponent of 
the decreasing ladder-height
process,
$\WH L^{-1}$ 
is the decreasing ladder-time
processes 
with
$\ell\doteq1/E[\WH L^{-1}(1)]$
(see Section~\ref{subsec:Setting}
for the definitions of 
$\WH\phi$ 
and
$\WH L^{-1}$),
$\Gamma(\cdot)$ 
denotes the gamma function
and the constant 
$C_\gamma$
is given by
\begin{equation}
\label{eq:cra_const}
C_\g \doteq \frac{\phi(0)}{\gamma \phi'(-\gamma)}.
\end{equation}
\end{Cor}

\subsection*{A brief description of the proofs and related literature}
The main result of this paper is the asymptotic independence
in Theorem~\ref{thm}.
Its proof, carried out 
in Section~\ref{sec:proof2},
is in three steps.
The first and in a certain sense the most
important step establishes the asymptotic behaviour of
the probability of an event involving the local time 
at zero of the reflected process. This event
contains precisely the paths sketched 
in Figure~\ref{fig:3_functionals}.
In the second step, a splitting property of an extended excursion process
of the reflected process,
introduced in the classical paper~\cite{GreenwoodPitman},
is applied to factorise the probabilities of certain events, related
to the ones involving the three functionals but with the running time
$t$ replaced by an independent exponential time $e(q)$.
The third and final step in the proof applies the
factorisations obtained in step two and the 
asymptotics from the first step to establish the 
stated asymptotic independence. This is achieved 
by studying (from first principles) the Laplace inversions 
of the probabilities 
arising in step two.

The law of the asymptotic overshoot, given by~\eqref{eq:zinf} of Theorem~\ref{thm}, 
is established in two steps in Section~\ref{sec:proof}. First, 
in Proposition~\ref{prop:aso} we extend the Cram\'er asymptotics
in~\cite{BertoinDoney} to the case of the two-sided exit. 
This step is based on the main result in~\cite{BertoinDoney}
and a renewal argument from~\cite{BertoinvanHarnSteutel}, applied in our setting under the 
Cram\'er measure. In the second step,
the law of 
$Z(x)$
is expressed in terms of the excursion measure of 
the reflected process. The limit law is then studied under
the excursion measure, using tools such as the result from step one, excursion
theory and the asymptotics from~\cite{DoneyMaller}.

The proof of Corollary~\ref{cor:4} is straightforward, once Theorem~\ref{thm}
has been established. Apart from our main result, it involves a number of classical facts 
from the fluctuation theory of L\'evy processes~\cite{Bertoin}. The details are given in 
Appendix~\ref{subsec:Cor}.

\section{Proof of Theorem~\ref{thm}}\label{sec:proof_Levy}

The asymptotic independence in Theorem~\ref{thm}
is a consequence of 
Proposition~\ref{prop:ind} proved in Section~\ref{sec:proof2}.
The formula for the law of the asymptotic overshoot 
follows from Lemma~\ref{lem:np} and Proposition~\ref{lem:exc}
established in Section~\ref{sec:proof}.
Section~\ref{subsec:Setting}
briefly defines the setting and notation of the proof.

\subsection{Setting and notation}
\label{subsec:Setting}
Let $(\Omega,\mc F,\{\mc F(t)\}_{t\ge 0}, P)$ be a filtered probability
space that carries a L\'{e}vy process
$X$
satisfying As.~\ref{A2}.
Here
$\Omega \doteq D(\mbb R)$ is
the Skorokhod space of real-valued functions that are
right-continuous on $\mbb R_+$ and have left-limits on $(0,\infty)$,
$X$
is the coordinate process,
$\{\mc F(t)\}_{t\ge 0}$
denotes the completed filtration generated by $X$,
which is right-continuous, and $\mc F$ is the completed
sigma-algebra generated by $\{X(t)\}_{t\ge 0}$.
For any $x\in\mbb R$ denote by $P_x$ the
probability measure on $(\Omega,\mc F)$ under which
$X-x$ is a L\'{e}vy process.
We refer to~\cite[Ch.~I]{Bertoin}
for further background on L\'{e}vy processes.

Let $L$ be a local time  at zero
of the reflected process $\WH Y=\{\WH Y(t)\}_{t\ge 0}$ of
the dual 
$\WH X \doteq -X$,
i.e.
$\WH Y(t) \doteq X^*(t)- X(t),$
where
$X^*(t)\doteq\sup_{0\leq s\leq t}X(s)$.
The ladder-time process
$L^{-1}~=~\{L^{-1}(t)\}_{t\ge 0}$
is equal to the right-continuous inverse of 
$L$.
The ladder-height process
$H=\{H(t)\}_{t\ge 0}$ is given by $H(t) \doteq 
X(L^{-1}(t))$
for all $t\ge 0$ with $L^{-1}(t)$ finite and by
$H(t)\doteq+\infty$
otherwise.
Let $\phi$ be the Laplace
exponent of $H$,
\begin{equation}\label{eq:phi}
\phi(\theta)\doteq-\log E[\te{-\theta H(1)}\I_{\{H(1)<\infty\}}],
\qquad\text{for any}\quad \theta\in\RR_+,
\end{equation}
where 
$\I_A$ denotes the indicator of a set $A$.
Analogously, define the local time 
$\WH L$
of $Y$
at zero,
the decreasing ladder-time and ladder-height subordinators
$\WH L^{-1}$
and
$\WH H$
with
$\WH \phi$
the Laplace exponent of 
$\WH H$.
See~\cite[Sec.~VI.1]{Bertoin}
for more details on ladder subordinators.
Note that the Cram\'{e}r assumption implies $E[X(1)]<0$,
making
$Y$
(resp. $\WH Y$)
a recurrent
(resp. transient) Markov process on
$\RR_+$. Hence
$\phi(0)>0$
and the stopping time $\tau(x)$
is a.s. finite for any
$x\in\RR_+$,
making
$H$
a killed subordinator
under 
$P$
and the overshoot
$Z(x)$
a $P$-almost surely defined random variable.

We now briefly review elements of It\^{o}'s excursion
theory that will be used in the proof. We refer to~\cite{GreenwoodPitman} 
and~\cite[Ch.~IV]{Bertoin} for a general treatment 
and further references.  Consider the Poisson point process 
of excursions away from zero associated to the strong Markov process
$Y$.
For each moment $t\in\mbb R_+$ of local time, let 
$\e(t)\in \mc E=\{\ve\in\Omega:\ve\ge 0\}$ 
denote the excursion at $t$: 
\begin{equation}\label{eq:e}
\e(t) \doteq
\begin{cases}
\le\{Y\le({s + \WH L^{-1}(t-)}\ri),\ s \in [0, \WH L^{-1}(t)-\WH L^{-1}(t-))\ri\}, & \text{if }
\WH L^{-1}(t-)<\WH L^{-1}(t),\\
\partial, & \text{otherwise},
\end{cases}
\end{equation}
where
$\partial\equiv 0$
is the null function,
$\WH L^{-1}({t-})\doteq \lim_{s\uparrow t}\WH L^{-1}(s)$
if $t>0$
and
$\WH L^{-1}(0-)=0$
otherwise.
Definition~\eqref{eq:e}
uses the fact
$\WH L(\infty)\doteq \lim_{s\to\infty}\WH L(s) =\infty$
$P$-a.s.,
which holds by 
the recurrence
of 
$Y$.
It\^{o}~\cite{Ito} proved that
$\e$ is a Poisson point
process under $P$.
Let
$n$
be the intensity (or excursion) measure on
$(\mc E, \mc G)$
of 
$\epsilon$, where $\mc G=\sigma(\epsilon(t), t\ge 0)$.
In Sections~\ref{sec:proof2} and~\ref{sec:proof},
for any Borel-measurable $F:\mc E \to \mbb R$ we
denote $n(F)=n(F(\ve)) \doteq \int_{\mc E}F\>\td n$. In this notation
the equality $n(A) = n(\I_A)$ holds for any $A\in\mathcal G$
and,
if 
$n(A)\in(0,\infty)$,
we denote $n(B|A) \doteq n(B\cap A)/n(A)$ for any $B\in\mc G$.

Let
$\zeta(\ve)\doteq\inf\{t> 0: \ve(t)=0\}$ denote the lifetime of an element
$\ve\in\mc E$.
(Note that
$\zeta(\e(t))$
is given in terms of
$\WH L^{-1}$
by
$\zeta(\e(t))= \WH L^{-1}(t) - \WH L^{-1}(t-)$
for any
$t\in\RR_+.$)
et $N^q$ be an independent standard
Poisson process with parameter $q$ and consider
the process $(X,N^q)$, which is defined on the probability space
$(\Omega\times\Omega, \mc F\otimes\mc F^N, P\times\mathbf P^N)$
where $\mc F^N$ is the completed filtration
generated by
$\{N^q(t)\}_{t\ge 0}$
and 
$\mbf P^N$
the probability law of $N^q$.
Let $\mathbb P \doteq P\times \mbf P^N$ be
the product measure and note that under $\PP$
the random variable
$T_{N^q}(1)$,
defined by
\begin{equation}
\label{eq:T_f}
T_f(1)\doteq\inf\{t\ge0: f(t)=1\}
\qquad\text{for any $f\in D(\mbb R)$,}
\end{equation}
is independent of $X$
and exponentially distributed with mean
$1/q$.
We associate to the L\'evy process
$(X,N^q)$ 
a two-dimensional point process
$(\e,\eta)=\{(\e(t),\eta(t))\}_{t\ge 0}$,
where
\begin{equation}
\label{eq:Def_eta}
\eta(t) \doteq
\begin{cases}
\le\{N^q({s + \WH L^{-1}(t-)})- N^q(\WH L^{-1}(t-)),\ s \in [0, \zeta(\e(t)))\ri\}, & \text{if }
\WH L^{-1}(t-)<\WH L^{-1}(t),\\
\partial, & \text{otherwise}.
\end{cases}
\end{equation}
Under $\PP$ the process $(\e,\eta)$ is a Poisson
point process with  values in $\mc E\times \mc E$. 
To the best of our knowledge, this construction first appeared in~\cite{GreenwoodPitman}.
We refer to
\cite[Ch. O.5]{Bertoin} for a treatment of Poisson point processes, 
the compensation formula and the properties of its characteristic measure.

\subsection{Asymptotic independence}
\label{sec:proof2}
The proof of the asymptotic independence in 
Proposition~\ref{prop:ind} below relies on the
following observations concerning the large time behaviour
of the local time $\WH L$.

\begin{Lemma}\label{lem:straddle}
The following statements hold true:
\begin{enumerate}
\item[(i)] The expectation of $\WH L^{-1}(1)$ satisfies $E[\WH L^{-1}(1)]\in(0,\infty)$. 
\item[(ii)] Recall from Corollary~\ref{cor:4} that $\ell=1/E[\WH L^{-1}(1)]$. 
Then for any 
$\delta\in(0,\ell/2)$
we have 
$$\limsup_{x\wedge t\to\infty}P(\WH L(\tau(x))\in t[\ell-\d,\ell+\d])\leq 
\frac{4}{\te{} \ell}\delta.$$ 
\item[(iii)] The following limit holds $P(\WH L(t) = \WH
L(\tau(x))) \longrightarrow 0$ as $x\wedge t\to\infty$;
\item[(iv)] For any $\delta_1,\delta_2\in[0,1/4)$ we have 
\begin{eqnarray}
\label{eq:Lemma_iv_ineq}
\limsup_{x\wedge t\to\infty}
P(\WH L(t(1-\delta_1)) \leq \WH L(\tau(x))\leq \WH L(t(1+\delta_2)))\leq
\frac{8}{\te{}}\max\{\delta_1,\delta_2\}.
\end{eqnarray}
For any fixed
$s\in(0,\infty)$ 
it holds $P(\WH L((t-s)\vee 0) \leq \WH L(\tau(x))< \WH L(t)) \longrightarrow 0$ 
as $x\wedge t\to\infty$.
\end{enumerate}
\end{Lemma}

\noindent{\it Remarks.} (1) Part~(iii) in Lemma~\ref{lem:straddle}
implies that,
as $x$ and $t$ tend to infinity,
the probability that the excursion straddling $t$
is the first excursion with height larger than $x$ tends to zero.
This fact can be viewed as an intuitive explanation for the
asymptotic independence of
$Z(\infty)$
and
$Y(\infty)$.
Part~(iv) 
of Lemma~\ref{lem:straddle} has analogous interpretation.\\
\noindent (2) The important role played by Lemma~\ref{lem:straddle} 
in the proof of the asymptotic independence in Theorem~\ref{thm}
lies in the fact that, the limits in parts~(iii) and~(vi) 
do not require the point
$(t,x)$
in
$(0,\infty)^2$
to tend to infinity along a specific trajectory  but only for 
its norm
$t\wedge x$
to increase beyond all bounds.\\
\noindent (3) In contrast to Lemma~\ref{lem:straddle}~(iv)
the inequality
$
\limsup_{x\wedge t\to\infty} P(\WH L(\tau(x))< \WH L(t) \leq \WH L(\tau(x+z)))>0
$
holds
for any fixed
$z>0$
(cf. remark following the statement of Lemma~\ref{lem:ind}).
To show this, recall 
$\WH L(t)/t\to \ell$
a.s.
as $t\uparrow\infty$
(see e.g. proof of Lemma~\ref{lem:straddle}~(iii) below)
and note that for any small $\delta>0$ we have
$P(\WH L(\tau(x))< \WH L(t) \leq \WH L(\tau(x+z)))\geq
P(\WH L(\tau(x))< t (\ell-\delta), \WH L(\tau(x+z))\geq t(\ell+\delta))+o(1)$.
Hence by Lemma~\ref{lem:np} and equality~\eqref{eq:DaneyMaller} we find
\begin{eqnarray*}
\lefteqn{P(\WH L(\tau(x))< t(\ell-\d), \WH L(\tau(x+z))\geq t(\ell+\d))}\\
&\geq&   P(\WH L(\tau(x+z))\geq t(\ell+\d)) - P(\WH L(\tau(x))\geq t(\ell-\d))\\
& = &  
\te{-t(\ell+\d)\, n(\rho(x+z)<\zeta)}-\te{-t(\ell-\delta)n(\rho(x)<\zeta)}
\rightarrow \te{-(\ell+\d)\lambda C_\gamma \WH \phi(\gamma)\te{-\gamma z}}-
\te{-(\ell-\d)\lambda C_\gamma \WH \phi(\gamma)}>0,
\end{eqnarray*}
where $x\wedge t\to\infty$
in such a way that
$t\te{-x\gamma}\to\lambda>0$.
Since
$z>0$,
the final inequality clearly holds for 
$\delta=0$
and hence by continuity for all 
$\delta>0$ 
sufficiently small.

\smallskip

\noindent {\it Proof of Lemma~\ref{lem:straddle}.} Part~(i)
of the lemma is known. For completeness a short proof, based on the 
Wiener-Hopf factorisation, is given in the Appendix.

(ii) For any $x,t\in(0,\infty)$, Lemma~\ref{lem:np} implies
$P(\WH L(\tau(x)) > t)
= \te{-t\, n(B(x))}$
for all
$t\geq0$,
where
$B(x)\doteq\{\rho(x)<\zeta\}$.
Therefore for any
$\delta\in(0,\ell/2)$
the following equality holds:
\begin{eqnarray*}
P(\WH L(\tau(x))\in t[\ell-\d,\ell+\d])
 =  \te{-t\,\ell\, n(B(x))}\left(\te{\delta\,t\, n(B(x))}-\te{-\delta\,t\, n(B(x))}\right).
\end{eqnarray*}
Lagrange's theorem implies that there exists
$\xi_{t,x}\in(-\delta,\delta)$
such that 
\begin{eqnarray*}
P(\WH L(\tau(x))\in t[\ell-\d,\ell+\d])
& = &  
2\delta t n(B(x))
\te{(\xi_{t,x}-\ell)t n(B(x))}\\
& \leq &
2\delta t n(B(x)) \te{-t n(B(x))\ell/2}\leq 
\delta 4/(\te{}\ell),
\end{eqnarray*}
where the inequality follows
from
$|\xi_{t,x}|<\ell/2$.
Since 
$t,x\in(0,\infty)$
were arbitrary,
this concludes the proof of part~(ii).

%

(iii) Since $\WH L^{-1}$ is a L\'{e}vy subordinator under $P$, the
strong law of large numbers
(see e.g.~\cite[p.92]{Bertoin}) implies that, as $t\to\infty$, the
ratio $t/\WH L^{-1}(t)$ tends to $\ell$ almost surely. Hence, for
any $\d\in(0,\ell/2)$,
\begin{equation}\label{lln}
P\le(\WH L(t)/t \in [\ell-\d,\ell+\d]\ri) = 1 + o(1), \qquad
\text{as $t\to\infty$}.
\end{equation}
Equation~\eqref{lln} 
yields the following
as
$t\wedge x\to\infty$:
\begin{eqnarray*}
P(\WH L(t) = \WH L(\tau(x))) &=& P(\WH L(t)
= \WH L(\tau(x)),
\WH L(t) \in t[\ell-\d,\ell+\d]) + o(1)\\
&\leq&
P(\WH L(\tau(x))\in t[\ell-\d,\ell+\d]) + o(1). 
\end{eqnarray*}
Hence part~(ii) yields
$\limsup_{x\wedge t\to\infty} P(\WH L(t) = \WH L(\tau(x)))\leq 
\delta 4/(\te{}\ell)$.
Since 
$\d\in(0,\ell/2)$
was arbitrary and probabilities are non-negative quantities,
the limit in part~(iii) follows.

(iv) Note that for any $\alpha\ge0$ the quotient $\WH
L(t\alpha)/t$ tends to $\ell\alpha$ $P$-a.s. as $t\to\infty$. For
any $\delta_1,\delta_2\in[0,1/4)$ we therefore find that the
probability of the event
$$A_{\delta_1,\delta_2}(t,x)=\{\WH L(t(1-\delta_1)) 
\leq \WH L(\tau(x))\leq \WH L(t(1+\delta_2))\} $$
satisfies the following as
$t\wedge x\to\infty$:
\begin{eqnarray}
\nonumber
P(A_{\delta_1,\delta_2}(t,x))
&=&
P(A_{\delta_1,\delta_2}(t,x),
\WH L(t(1-\delta_1)),
\WH L(t(1+\delta_2))\in t[\ell(1-\d),\ell(1+\d)]) + o(1)\\
\label{eq:Part_4_proof}
&\leq& P(\WH L(\tau(x))\in t[\ell(1-\d),\ell(1+\d)]) + o(1),
\end{eqnarray}
for any
$\delta\in(2\max\{\delta_1,\delta_2\},1/2)$.
Since 
$0<\delta \ell<\ell/2$,
part~(ii) of the lemma and inequality~\eqref{eq:Part_4_proof}
imply that 
$\limsup_{x\wedge t\to\infty} P(A_{\delta_1,\delta_2}(t,x))\leq 
\delta 4/\te{}$.
Therefore the first inequality in part~(iv) 
is satisfied.
The second limit in part~(iv) follows by noting that,
for any
$s\in\RR_+$
and
$\delta_1\in(0,1/4)$,
the inclusion
$ \{\WH L((t-s)\vee 0) \leq \WH L(\tau(x))< \WH L(t)\}\subset A_{\delta_1,0}(t,x)$
holds for all
$(t,x)$ with large
$t\wedge x$.
Hence by~\eqref{eq:Lemma_iv_ineq} 
the following holds
$$ \limsup_{t\wedge x\to\infty}
P(\WH L((t-s)\vee 0) \leq \WH L(\tau(x))< \WH L(t))\leq 
\delta_1 8/\te{}.$$
Since $\delta_1$
can be chosen arbitrarily small, this
proves part~(iv) 
and hence the lemma.~\exit

\smallskip

Recall that the random variable $e(q) \doteq T_{N^q}(1)$, 
where
$N^q$ 
is defined immediately above~\eqref{eq:T_f}
and
$T_{N^q}(1)$ 
(for any path of 
$N^q$)
is given in~\eqref{eq:T_f},
is exponentially distributed with mean $1/q$ 
and independent 
of $X$. We 
establish the following lemma:

\begin{Lemma}\label{lem:ind} For any $q,x \in (0,\infty)$,
$z\in[-x,0]$,
$y\in[0,x]$
and Borel sets $A,B,C \in \mc B(\mbb R_+)$
define the quantities
\begin{equation*}
\pi_1(q,A)\doteq\PP(Y({e(q)})\in A),\quad \pi_2(x,B)\doteq P(Z(x)\in B),\quad
\pi_3(q,x+z)\doteq\PP(\WH L({\tau(x+z)}) < \WH L({e(q)}))
\end{equation*}
and
\begin{eqnarray*}
r(y,x,q,C,B) & \doteq & \pi_2(x,B)\PP(Y(e(q))\in C, \WH L(\tau(y))=\WH L(e(q)))\\
& - & \PP(Y(e(q))\in C,Z(x)\in B, \WH L(\tau(y))=\WH L(e(q))),\\
R(q,x) &\doteq& \pi_1(q,A)\> r(x,x,q,\RR_+,B).
\end{eqnarray*}
Then the following equalities hold:
\begin{eqnarray}
&&\label{eq:idl1}\PP(Y({e(q)})\in A, Z(x)\in B) =
\pi_1(q,A)\pi_2(x,B) + R_1(q,x),\\
&&\label{eq:idl2}
\PP(Y({e(q)})\in A, Z(x)\in B,\WH L({\tau(x+z)}) < \WH L({e(q)}))
= \pi_1(q,A)\pi_2(x,B)\pi_3(q,x+z) + R_2(q,x,z),
\end{eqnarray}
where
\begin{eqnarray*}
R_1(q,x) &\doteq& R(q,x) - r(x,x,q,A,B),\\
R_2(q,x,z) & \doteq & R_1(q,x)+r(x+z,x,q,A,B). 
\end{eqnarray*}
\end{Lemma}


\noindent{\it Remark.} The proof of the asymptotic independence 
of the triplet 
$(Y(t),Z(x+y),M(t,x))$
in Proposition~\ref{prop:ind}
is based on~\eqref{eq:idl2} and the fact that 
$R_2(q,x,z)$ is a linear combination of the probabilities
of events, each of which is contained in an event of
the form 
$\{\WH L(e(q))=\WH L(\tau(x))\}$,
the probability of which tends to zero as $t\wedge x \to \infty$ 
(cf.  Lemma~\ref{lem:straddle}(iii)).
It is important to note that the equality in~\eqref{eq:idl2}
cannot be extended to the case
$z>0$, since the random variables
$\I_{\{ Z(x)\in B \}}$
and
$\I_{\{\WH L(\tau(x+z)) < \WH L(e(q))\}}$
are clearly functions of the same excursion 
on the event
$\{\WH L({\tau(x+z)}) = \WH L(\tau(x))\}$
consisting of the paths of 
$Y$
that cross the levels 
$x$
and
$x+z$
for the first time during the same excursion.
In particular
$\{Z(x)>z\}\subset\{\WH L({\tau(x+z)}) = \WH L(\tau(x))\}$
and hence
$P(\WH L({\tau(x+z)}) = \WH L(\tau(x)))$
is in the limit as 
$x\to\infty$
bounded below by the strictly positive probability 
of 
$\{Z(\infty)>z\}$.
This observation invalidates the proof of 
Lemma~\ref{lem:ind}
if 
$z>0$.
Furthermore, it is not difficult to see that
in general, for 
$z>0$,
the events
$\{M(t,x)\leq z\}$
and 
$\{Z(x)\in B\}$
are not asymptotically independent as
$t\wedge x\to\infty$.

\smallskip

\noindent {\it Proof of Lemma \ref{lem:ind}.}
Define the set
$A' \doteq \{(\ve,\mu)\in \mc E\times \mc E:T_\mu(1) < \zeta(\ve)\}$
and let
$H_{A'}=\inf\{t\ge0: (\epsilon(t),\eta(t))\in A'\}$
be the first entry of the Poisson point process
$(\epsilon,\eta)$,
defined in~\eqref{eq:e} and~\eqref{eq:Def_eta},
into the set 
$A'$. For any set $A\in\mc G$, denote by $\e^A$ the Poisson point process $\e$
killed upon its first entrance into $A$, i.e. $\e^A(t)$ is equal to $\e(t)$
if $t<H_A$ and to $\partial$ otherwise, and define $\eta^A$ analogously.

The
definitions of the Poisson point process
$(\epsilon,\eta)$
in~\eqref{eq:e} and~\eqref{eq:Def_eta}
and that of 
$T_f(1)$
in~\eqref{eq:T_f}
imply 
$e(q)= T_{N^q}(1) = \WH L^{-1}(H_{A'}-)+T_{\eta(H_{A'})}(1)<\WH L^{-1}(H_{A'})$,
which yields 
$Y(e(q)) = \epsilon(H_{A'})(T_{\eta(H_{A'})}(1))$
and
$\WH L(e(q))=H_{A'}.$ 
Hence
$\{Y({e(q)})\in A\}\in\s(\e(H_{A'}), \eta(H_{A'}))$
and the event 
$\{\WH L({\tau(x)}) < \WH L({e(q)}), Z(x)\in B\}$
is measurable with respect to the sigma-algebra
$\s(\e^{A'},\eta^{A'})$.
Therefore the two events are independent by the
splitting property (see~\cite[Sec~O.5, Prop.~O.2]{Bertoin})
of the Poisson point process
$(\e,\eta)$
(note for example that the latter event would not be measurable with respect to
$\s(\e^{A'},\eta^{A'})$
if the strict inequality was replaced by
``$\leq$'').
Hence 
the following equality
holds:
\begin{eqnarray}\label{eq:LT1}
\lefteqn{\PP(\WH L({\tau(x)}) < \WH L({e(q)}), Z(x)\in B, Y({e(q)})\in A)}\\  
& =  &
\pi_1(q,A) \PP(\WH L({\tau(x)}) < \WH L({e(q)}), Z(x)\in B).
\nonumber
\end{eqnarray}
An analogous argument
based on the splitting property
of the Poisson point process
$(\epsilon,\eta)$
implies
that  the events
$\{Z(x)\in B\}$ 
and 
$\{\WH L({\tau(x)}) \leq \WH L({e(q)})\}$ are independent. Indeed,
let
$B'\doteq\{(\ve,\mu)\in \mc E\times\mc E:\rho(x,\ve)<\zeta(\ve)\}$
and note that
$\WH L(\tau(x))=H_{B'}$,
which implies that 
$\{\WH L(\tau(x)) > \WH L(e(q))\}=\{H_{B'}>H_{A'}\}$
and hence 
$\{\WH L(\tau(x)) \leq \WH L(e(q))\}\in\s(\e^{B'},\eta^{B'})$.
Furthermore 
$
\{Z(x)\in B\}  
 =  
\left\{
\epsilon(H_{B'})(\rho(x,\epsilon(H_{B'})))-x\in B
\right\}\in\s(\e(H_{B'}))
$
and hence the
splitting property~\cite[Sec~O.5, Prop.~O.2]{Bertoin}
of the Poisson point process
$(\e,\eta)$
implies
\begin{eqnarray}\label{eq:LT2}
\PP(\WH L({\tau(x)}) \leq \WH L({e(q)}), Z(x)\in B)
& = &
\PP(\WH L({\tau(x)}) \leq \WH L({e(q)}))\>P(Z(x)\in B).
\end{eqnarray}
The identities in~\eqref{eq:LT1} and~\eqref{eq:LT2}
imply the equality in~\eqref{eq:idl2}
in the case 
$z=0$.



To prove~\eqref{eq:idl1} we first observe that
the equality 
\begin{equation}\label{eq:LT3}
\PP(\WH L({\tau(x)}) > \WH L({e(q)}),  Y({e(q)})\in A, Z(x) \in B)
 =  \PP(\WH L({\tau(x)}) > \WH L({e(q)}),  Y({e(q)})\in A)\pi_2(x,B)\\ 
\end{equation}
follows from the splitting property of the Poisson point process
$(\epsilon,\eta)$
at
$H_{B'}$.
The identity
\begin{eqnarray}
\label{eq:LT4}
\PP(\WH L({\tau(x)}) > \WH L({e(q)}),  Y({e(q)})\in A) & = &
\PP(\WH L({\tau(x)}) \ge \WH L({e(q)})) \PP( Y({e(q)})\in A)
\\
& - &
\PP(\WH L({\tau(x)}) = \WH L({e(q)}),  Y({e(q)})\in A) 
\nonumber
\end{eqnarray}
follows 
by~\eqref{eq:LT1}
for
$B=\RR_+$
and
$z=0$.
The equality in~\eqref{eq:idl2} for 
$z=0$,
which was proved above,
and identities~\eqref{eq:LT3} and~\eqref{eq:LT4}
imply the equality in~\eqref{eq:idl1}. This concludes the proof of~\eqref{eq:idl1}. 

For any
$z\in[-x,0)$
it holds
$\WH L(\tau(x+z))\leq\WH L(\tau(x))$
and hence,
by the splitting property of the Poisson point process
$(\epsilon,\eta)$
at
$H_{B'}$,
we find
\begin{eqnarray}\label{eq:LT3_z_negative}
\lefteqn{\PP(\WH L({\tau(x+z)}) > \WH L({e(q)}),  Y({e(q)})\in A, Z(x) \in B)}\\
& = &  \PP(\WH L({\tau(x+z)}) > \WH L({e(q)}),  Y({e(q)})\in A)\pi_2(x,B).
\nonumber
\end{eqnarray}
Furthermore the splitting property
of 
$(\epsilon,\eta)$
applied one more time at
$H_{A'}$
yields
\begin{eqnarray}\label{eq:LT1_z_negative}
\PP(\WH L({\tau(x+z)}) < \WH L({e(q)}), Y({e(q)})\in A)
& =  &
\pi_1(q,A) \PP(\WH L({\tau(x+z)}) < \WH L({e(q)})).
\end{eqnarray} 
The following elementary equality
\begin{eqnarray*}
\lefteqn{\PP(\WH L({\tau(x+z)}) < \WH L({e(q)}),  Y({e(q)})\in A, Z(x) \in B) 
 =  \PP(Y({e(q)})\in A, Z(x) \in B)}\\ 
& - & 
\PP(\WH L({\tau(x+z)}) > \WH L({e(q)}),  Y({e(q)})\in A, Z(x) \in B) \\
& - & 
\PP(\WH L({\tau(x+z)}) = \WH L({e(q)}),  Y({e(q)})\in A, Z(x) \in B)
\end{eqnarray*} 
and
identities~\eqref{eq:idl1},~\eqref{eq:LT3_z_negative},~\eqref{eq:LT1_z_negative}
imply~\eqref{eq:idl2}
holds 
in the case 
$z\in[-x,0)$.
This concludes the proof of the lemma.~\exit

\smallskip

Before proving the asymptotic independence
of
$Y(t),Z(x+y)$
and
$M(t,x)$
stated in 
Proposition~\ref{prop:ind} below,
we need to establish the asymptotic behaviour of certain convolutions
that will arise in the proof 
of
Proposition~\ref{prop:ind}.
Let $T(x)$ and $\WH T(x)$ denote the
first-passage times of $X$ into
the intervals
$(x,\infty)$ and $(-\infty,-x)$
respectively for any
$x\geq0$,
\begin{equation}\label{eq:TT}
T(x)\,\doteq\, \inf\{t\ge 0: X(t)\in(x,\infty)\},\qquad
\WH T(x)\,\doteq\, \inf\{t\ge0: X(t)\in(-\infty,-x)\}.
\end{equation}

\begin{Lemma}\label{lem:InversLaplaceProduct}
Let
$a\in [0,\infty)$
and recall that
$T(a)$
is the first-passage time of
$X$
over the level
$a$
defined in
\eqref{eq:TT}.
Then the following equality holds:
\begin{equation}
\label{eq:Limit_of_Int_1}
\int_{[0,t]}P(\WH L({\tau(y)}) = \WH L(t-s)) \,P(T(a)\in\td s)  =  o(1)
\end{equation}
as $y\wedge t\to \infty$.
Furthermore, we have
\begin{equation}
\label{eq:Limit_of_Int_2}
\int_{[0,t]}P(\WH L({\tau(y)}) < \WH L(t-s))\,P(T(a)\in\td s)
=P(\WH L({\tau(y)}) < \WH L(t)) \, P(T(a)\leq t) + o(1)
\end{equation}
as $y\wedge t\to \infty$.
\end{Lemma}

\noindent {\it Proof of Lemma~\ref{lem:InversLaplaceProduct}.}
The proof of this lemma is based on Lemma~\ref{lem:straddle}.
Note that for fixed
$t,y\in(0,\infty)$,
the integral in~\eqref{eq:Limit_of_Int_1}
can be expressed as an integral over
$\RR_+$
(with respect to the measure
$P(T(a)\in\td s)$)
of the integrand
$s\mapsto I_{[0,t]}(s)P(\WH L({\tau(y)}) = \WH L(t-s))$.
Lemma~\ref{lem:straddle}~(iii) implies
that for any fixed
$s\in\RR_+$
the integrand tends to zero
as $y\wedge t\to \infty$.
Therefore~\eqref{eq:Limit_of_Int_1}
follows as a consequence of the dominated convergence theorem,
since the integrands are uniformly bounded by
one and the measure is finite.

To prove equality~\eqref{eq:Limit_of_Int_2},
first note that it 
is equivalent to the statement
\begin{equation}
\label{eq:Step_Between}
\int_{[0,t]}\left(P(\WH L({\tau(y)}) < \WH L(t))-P(\WH L({\tau(y)})<\WH L(t-s))\right)P(T(a)\in\td
s)= o(1)
\end{equation}
as $y\wedge t\to \infty$.
Since the local time
$\WH L$
is non-decreasing, the integrand in~\eqref{eq:Step_Between}
can be expressed as
\begin{equation}
\label{eq:Step_Between_2}
P(\WH L({\tau(y)}) < \WH L(t))- P(\WH L({\tau(y)})<\WH L(t-s)) =
P(\WH L(t-s)\leq \WH L({\tau(y)}) < \WH L(t)).
\end{equation}
Equality~\eqref{eq:Step_Between_2},
Lemma~\ref{lem:straddle}~(iv)
and the dominated convergence theorem imply that~\eqref{eq:Step_Between},
and hence~\eqref{eq:Limit_of_Int_2}, holds.
This completes the proof of the lemma.
\exit

\begin{Prop}\label{prop:ind}Let
$A=(a,\infty)$
for some
$a\in\RR_+$,
$B \in \mc B(\mbb R_+)$
and
$C=(-\infty,z]$
for 
$z\in\RR$.
Then the following holds as 
$t\wedge y\wedge(x-y)\to\infty$:
\begin{eqnarray}
\label{eq:id1}
&& P(Y(t)\in A, Z(x)\in B) =
P(Y(t)\in A)P(Z(x)\in B) + o(1),
\\
\label{eq:id2}
&& P(Y(t)\in A, Z(x)\in B, M(t,y)\in C)
=  P(Y(t)\in A)P(Z(x)\in B)P(M(t,y)\in C) + o(1).
\end{eqnarray}
\end{Prop}

\noindent {\it Proof of Proposition \ref{prop:ind}.}
We first prove equality~\eqref{eq:id1}. 
Note that
$t\wedge y\wedge(x-y)\to\infty$
in particular implies
$t\wedge x\to\infty$
and
$t\wedge y\to\infty$.
By~\eqref{eq:idl1}
we have
\begin{eqnarray*}
P(Y(t)\in A, Z(x)\in B) & =  &
\mathcal L^{-1}\left(q\mapsto \frac{1}{q}\PP(Y(e(q))\in A, Z(x)\in B)\right)(t)\\
& = & P(Y(t)\in A)P(Z(x)\in B) +
\mathcal L^{-1}\left(q\mapsto \frac{1}{q}R_1(q,x)\right)(t),
\end{eqnarray*}
where
$\mathcal L^{-1}$
denotes the inverse Laplace transform and
$R_1(q,x)$
is defined in Lemma~\ref{lem:ind}.
Furthermore, 
\begin{eqnarray}
\nonumber
\mathcal L^{-1}\left(q\mapsto \frac{1}{q}R_1(q,x)\right)(t)
& = &
\mathcal L^{-1}\left(q\mapsto \frac{1}{q}R(q,x)\right)(t) -
P(Y(t)\in A,Z(x)\in B,\WH L({\tau(x)})=\WH L(t))\\
\nonumber
& + &
P(Z(x)\in B) P(Y(t)\in A,\WH L({\tau(x)})=\WH L(t))\\
\label{eq:R_1_is_R_asym}
& = &
\mathcal L^{-1}\left(q\mapsto \frac{1}{q}R(q,x)\right)(t) + o(1),
\qquad\text{as}\quad
t\wedge x\to \infty,
\end{eqnarray}
where the second equality holds
by Lemma~\ref{lem:straddle}~(iii).

To prove~\eqref{eq:id1}
we therefore need to establish the equality 
\begin{equation}
\label{eq:R_1_is_o_1} \mathcal L^{-1}\left(q\mapsto \frac{1}{q}R(q,x)\right)(t) = o(1)
\qquad\text{as}\quad
t\wedge x\to \infty.
\end{equation}
Since 
for every
$t$,
$Y(t)$
has the same law as 
$X^*(t)=\sup_{0\leq s\leq t}X(s)$ ,
$P(\Delta X(t)=0)=1$
for all 
$t>0$,
where
$\Delta X(t)\doteq X(t)-X(t-)$,
and
$\{X^*(t)>a,\Delta X(t)=0\}=\{T(a)<t,\Delta X(t)=0\}$,
the following equalities hold:
\begin{eqnarray}
\label{eq:Lap_of_Y_e_q}
\PP(Y(e(q))\in A)  =  \PP(X^*(e(q))>a)=\PP(T(a)<e(q)) 
 =  \int_{[0,\infty)}\te{-qt}P(T(a)\in\td t),
\end{eqnarray}
where as before
$e(q)$ 
is an exponential time
with mean
$1/q$ that is independent of $X$
and
$T(a)$
is defined in~\eqref{eq:TT}.
Since 
$q\mapsto \PP(Y(e(q))\in A)$
is by~\eqref{eq:Lap_of_Y_e_q}
the Laplace transform of the positive measure 
$P(T(a)\in\td t)$
on
$\RR_+$,
the following holds:
\begin{eqnarray*}
\lefteqn{\mathcal L^{-1}\left(q\mapsto \PP(Y(e(q))\in A) \frac{1}{q}
\PP(\WH L({\tau(x)})=\WH L(e(q))) \right)(t)}\\
& = &
\int_{[0,t]}P(\WH L({\tau(x)})=\WH L(t-s)) \,P(T(a)\in\td s)
 =   o(1)
 \qquad\text{as}\quad x\wedge t\to \infty,
\end{eqnarray*}
where the final equality follows 
by~\eqref{eq:Limit_of_Int_1} in Lemma~\ref{lem:InversLaplaceProduct}.
An analogous argument shows that
$$\mathcal L^{-1}\left(q\mapsto \PP(Y(e(q))\in A) \frac{1}{q}
\PP(Z(x)\in B, \WH L({\tau(x)})=\WH L(e(q))) \right)(t)=o(1)
\qquad\text{as}\quad x\wedge t\to \infty.$$
The definition of
$R(q,x)$
in Lemma~\ref{lem:ind}
and the two equalities above
imply~\eqref{eq:R_1_is_o_1}
and hence~\eqref{eq:id1}.

The proof of~\eqref{eq:id2}
is based on 
equality~\eqref{eq:Intermediate_Step} below,
which we now establish.
Since by assumption
$t\wedge y\wedge(x-y)\to\infty$,
for large values of 
$x$
and
$y$
we have 
$0\leq y\leq x$.
The definition of
$r(y,x,q,A,B)$
in Lemma~\ref{lem:ind} implies
\begin{eqnarray}
\nonumber
\mathcal L^{-1}\left(q\mapsto \frac{1}{q}r(y,x,q,A,B)\right)(t)
& = &  P(Z(x)\in B) P(Y(t)\in A,\WH L({\tau(y)})=\WH L(t))\\ 
\nonumber
& - &
P(Y(t)\in A,Z(x)\in B,\WH L({\tau(y)})=\WH L(t))\\
& = & o(1)
\qquad\text{as}\quad
t\wedge y\wedge(x-y)\to\infty,
\label{eq:r_asymp}
\end{eqnarray}
where the final equality follows from
Lemma~\ref{lem:straddle}~(iii).
The identity in~\eqref{eq:idl2} together with~\eqref{eq:r_asymp}, 
identities~\eqref{eq:R_1_is_R_asym},~\eqref{eq:R_1_is_o_1}
and~\eqref{eq:Lap_of_Y_e_q}
and equality~\eqref{eq:Limit_of_Int_2} in Lemma~\ref{lem:InversLaplaceProduct}
imply:
\begin{eqnarray*}
\lefteqn{P(Y(t)\in A, Z(x)\in B,\WH L({\tau(y)}) < \WH L({t}))}\\
&= & \mathcal L^{-1}\left(q\mapsto P(Z(x)\in B) \PP(Y(e(q))\in A) \frac{1}{q}
\PP(\WH L({\tau(y)})<\WH L(e(q))) +
\frac{1}{q}R(q,x)\right)(t)  + o(1)
\\
&=& P(Z(x)\in B)P(T(a)\leq t)P(\WH L({\tau(y)}) < \WH L({t})) + o(1)
\qquad\text{as}\quad 
t\wedge y\wedge(x-y)\to\infty.
\end{eqnarray*}
The process
$X$
drifts to 
$-\infty$
as
$t\to\infty$
by As.~\ref{A2},
which implies
$\lim_{t\to\infty}P(T(a)= t)=0$.
Thus, we have the following equality for the set $A=(a,\infty)$:
$$P(T(a)\leq t)=P(T(a)< t)+P(T(a)= t)=P(Y(t)\in A)+o(1)
\qquad \text{as $t\to\infty$.}
$$
As a consequence the following asymptotic independence holds:
\begin{eqnarray}
\nonumber
\lefteqn{P(Y(t)\in A, Z(x)\in B,\WH L({\tau(y)}) < \WH L({t}))}\\
&=& P(Y(t)\in A)P(Z(x)\in B)P(\WH L({\tau(y)}) < \WH L({t})) + o(1)
\qquad\text{as}\quad 
t\wedge y\wedge(x-y)\to\infty.
\label{eq:Intermediate_Step}
\end{eqnarray}

Recall that 
$C=(-\infty,z]$
for an arbitrary fixed
$z\in\RR$.
In order to prove equality~\eqref{eq:id2}
note that the following inclusions hold for any
$y\in\RR_+$:
\begin{eqnarray*}
\{M(t,y) \in C\}=\{Y^*(t)\leq y+z\} & \subset &  \{\WH L(t) \leq \WH L({\tau((y+z)^+)})\}\quad\text{and}
\\
\{\WH L(t) \leq \WH L({\tau((y+z)^+)})\} \cap \{M(t,y) \notin C\}
&  \subset &
\{\WH L({\tau((y+ z)^+)}) = \WH L(t)\}
\end{eqnarray*}
(recall that $\tau(x)$
is defined for
$x\in\RR_+$).
These inclusions, together with
Lemma~\ref{lem:straddle}~(iii),
imply that the following equality holds
for any family of events
$E(t,x)\in\mathcal F$,
$t,x\in\RR_+$,
as
$t\wedge y\wedge(x-y)\to\infty$:
\begin{eqnarray}
\label{eq;E_equality_F}
P\left(E(t,x),\WH L(t) \leq \WH L(\tau((y+z)^+))\right)
& = &
P\left(E(t,x),M(t,y)\in C\right) + o(1).
\end{eqnarray}
Since
$t\wedge y\wedge(x-y)\to\infty$,
for the fixed
$z\in\RR_+$
the inequalities 
$0\leq y+z\leq x$
hold for all large $y$ and $x$.
In particular~\eqref{eq:Intermediate_Step},
applied to the complement
$\{\WH L({\tau(y+z)}) < \WH L({t})\}^c= \{\WH L({\tau(y+z)}) \geq \WH L({t})\}$,
and~\eqref{eq;E_equality_F}
yield the following equalities
\begin{eqnarray*}
P(Y(t)\in A, Z(x)\in B,M(t,y)\in C) & = &
P(Y(t)\in A, Z(x)\in B,\WH L({t})\leq  \WH L({\tau(y+z)}))
+ o(1)\\ 
& = &
P(Y(t)\in A) P(Z(x)\in B)P(\WH L({t})\leq  \WH L({\tau(y+z)})) + o(1)\\
& = & P(Y(t)\in A) P(Z(x)\in B)P(M(t,y)\in C)+ o(1)
\end{eqnarray*}
as
$t\wedge y\wedge(x-y)\to\infty$.
This concludes the proof of~\eqref{eq:id2}.~\exit

\subsection{Limiting overshoot}
\label{sec:proof}
In this section we prove the formula in~\eqref{eq:zinf}
of Theorem~\ref{thm}, which characterises the law of the
limiting overshoot
$Z(\infty)$.
This is achieved in two steps. 
We first establish Cram\'er's asymptotics for the exit probabilities of 
$X$
from a finite interval. 
In the second step we describe the distribution of the overshoot 
$Z(x)$,
defined in~\eqref{eq:Levy_Overshoot},
in terms of the excursion measure $n$
(see Sections~\ref{sec:proof2} for the definition of $n$)
and apply the result from step one to find the relevant asymptotics 
under the excursion measure, which in turn yield 
the Laplace transform of the limiting law
$Z(\infty)$.

For any $x\in\RR_+$,  
recall that $T(x)$ is given in~\eqref{eq:TT} 
and define the overshoot
$$K(x)\doteq X({T(x)}) - x \qquad\text{ on the event $\{T(x) < \infty\}$.}$$
Denote by
$f(x)\simeq g(x)$
as
$x\uparrow\infty$
the functions
$f,g:\mbb R_+\to(0,\infty)$
satisfying
$\lim_{x\uparrow\infty}\frac{f(x)}{g(x)}=1$.

\begin{Prop}
\label{prop:aso} (i) {\bf (Asymptotic
two-sided exit probability)} For any $z>0$ we have
\begin{equation}
\label{eq:Prob_Asym_LEvy}
P(T(x) < \WH T(z)) \simeq C_\gamma\te{-\gamma x}\le(1 -
E\le[\te{\gamma X({\WH T(z)})}
\ri]\ri)
 \qquad\text{as}\quad x\to\infty,
\end{equation}
where  the constant $C_\gamma$ is given in~\eqref{eq:cra_const}
and $\WH T(z)$ in~\eqref{eq:TT}.

\noindent (ii) {\bf (Asymptotic overshoot)} Let $u\in\mathbb R_+$ and fix
$z>0$. Then we have
as
$x\to \infty$:
\begin{equation}
\label{eq:Most_Important} E\left[\te{-u K(x)}\I_{\{T(x)<\WH
T(z)\}}\right] \simeq C(u)\te{-\gamma x} \le( 1 -
E\left[\te{\gamma X(\WH T(z))}
\right]\ri), \quad \text{with}\quad
C(u) \doteq \frac{\gamma}{\gamma + u}\cdot \frac{\phi(u)}{\phi(0)}\cdot C_\gamma 
\end{equation}
and $C_\gamma$ 
in~\eqref{eq:cra_const}. 
\end{Prop}


\noindent{\it Remarks.} (i) 
Let
$P^{(\gamma)}_x$
be the Cram\'{e}r measure
on $(\Omega,\mc F)$.
Its restriction to 
$\mc F(t)$ is given by
\begin{equation*}
P_x^{(\gamma)}(A) \doteq E_x[\te{\gamma (X(t)-x)}\I_A],\qquad A\in\mc F(t),\qquad t\in\mbb R_+.
\end{equation*}
Here
$E_x$
is the expectation under
$P_x$
and $\I_A$ is the indicator of 
$A$.
Under As.~\ref{A2} it follows that
$P^{(\gamma)}_x$
is a probability measure and
$X-x$
is a L\'{e}vy process under
$P_x^{(\gamma)}$
with
$E^{(\gamma)}_x[X(1)-x]\in(0,\infty)$.

\noindent (ii) Since the overshoot of $X$ is the same as that 
of its ladder process, the weak limit under $P^{(\gamma)}$ of $K(x)$ as $x\to\infty$,
needed in the proof of Proposition~\ref{prop:aso}
is be derived from \cite[Thm. 1]{BertoinvanHarnSteutel}.

\noindent (iii) Note that the random variable $X({\WH T(z)})$ under the
expectation in~\eqref{eq:Prob_Asym_LEvy} is well-defined $P$-a.s.,
since As.~\ref{A2} implies that the L\'evy process $X$
drifts to $-\infty$ $P$-a.s.

\noindent (iv) The proof of the Proposition~\ref{prop:aso} 
is based on two ingredients: the Cram\'er estimate for L\'evy processes~\cite{BertoinDoney}
and the fact that the overshoot 
$K(x)$
has a weak limit under 
$P^{(\gamma)}$ 
follows from~\cite[Thm.~1]{BertoinvanHarnSteutel}.
The details of the proof are 
given in Appendix~\ref{subsec:Two_sided_Overshoot_LEvy}.


\smallskip

Let 
$x>0$
and denote 
by
$\rho(x,\ve)$
the first time that an excursion
$\ve\in\mc E$
enters the interval
$(x,\infty)$:
\begin{equation}\label{eq:rho}
\rho(x,\ve) \doteq \inf\{s\ge0: \ve(s) > x\}.
\end{equation}
For brevity we sometimes write $\rho(x)$ instead of $\rho(x,\ve)$.
Since the expectation
$E^{(\gamma)}[X_1]$
is strictly positive, 
under 
$P^{(\gamma)}$
the reflected process $Y$ is transient 
and
$\WH L(\infty)$
is an exponentially distributed 
random variable, independent of the killed subordinator 
$\{(\WH L^{-1}(t),\WH H(t))\}_{t\in[0,\WH L(\infty))}$.
As a consequence, the excursion process $\e'=\{\e'(t)\}_{t\ge 0}$,
defined by the formula in~\eqref{eq:e}
for $t<\WH L({\infty})$ and
by $\e'(t)\doteq\partial$ otherwise, is under
$P^{(\gamma)}$
a Poisson point process
killed at an independent exponential time
with mean $E^{(\gamma)}[\WH L({\infty})]$.
Put differently,
$\e'$
is given by~\eqref{eq:e}
up to the first time it hits the set
$\{\varepsilon\in\mathcal E:\zeta(\ve)=\infty\}$.
In the rest of the paper we will
denote by
$n^{(\gamma)}$
the excursion measure
under $P^{(\gamma)}$
of the killed Poisson point process
$\e'$.

\begin{Lemma}\label{lem:np}
For any
$x>0$
the random variable
$\WH L({\tau(x)})$
is exponentially distributed
under
$P$
(resp.
$P^{(\gamma)}$)
with parameter
$n(\rho(x)<\zeta)$
(resp. $n^{(\gamma)}(\rho(x)<\zeta)$)
and the following equality holds:
\begin{eqnarray}
\nonumber
P( Z(x) > y) &=& n(\ve(\rho(x,\ve)) - x > y|\rho(x)<\zeta)\qquad\text{for any}\quad
y\in\RR_+.
\end{eqnarray}
\end{Lemma}

\proof
The definitions of the Poisson point process $\e$ in~\eqref{eq:e} and
the first-passage time
$\rho(x,\ve)$
in~\eqref{eq:rho}
imply the equality
$\WH L({\tau(x)}) = H_A \doteq \inf\{t\ge 0: \epsilon(t)\in A\}$
where
$A \doteq \{\ve\in\mc E: \rho(x,\ve)< \zeta(\ve)\}$.
The first
statement in the lemma follows since $H_A$ is
exponentially distribution with parameter $n(A)$ 
(e.g.~\cite[Sec.~O.5, Prop.~O.2]{Bertoin}).
The second statement is a consequence of 
the fact that $\e(H_A)$ follows an $n$-uniform distribution
(i.e. $P(\e(H_A) \in B)=n(B|A)$ for any $B\in\mc G$, see 
e.g.~\cite[Sec.~O.5, Prop.~O.2]{Bertoin}),
taking $B$ to be equal to
$\{\ve\in \mc E: \rho(x,\ve) < \zeta(\ve), \ve(\rho(x,\ve)) - x > y\}$.
\exit

Conversely, one may also express $n$ as a ratio of expectations
under the measure $P$. To derive such a representation, 
for any
$x>0$,
define the random variable $K_F(x)$ by
\begin{equation}
\label{eq:Def_K_F}
K_F(x) \doteq \sum_{g} F(\e_g)\I_{\{g<\tau(x)\}},
\end{equation}
where the sum runs over all left-end points $g$ of
excursion intervals,
$\e_g\doteq \e (\WH L(g))$,
and
$F:\mc E \to \mbb R$
is Borel-measurable and non-negative
(note that
$F\equiv1$
implies
$K_F(x)\equiv1$
$P$-
and 
$P^{(\gamma)}$-almost surely).

\begin{Lemma}\label{lem:npp} (i) 
Define
$\WH{\mc V}({x}) \doteq E\le[\WH L(\tau(x))\ri]$
and
$\WH{\mc V}^{(\gamma)}({x}) \doteq E^{(\gamma)}\le[\WH L(\tau(x))\ri]$.
Then the following hold:
\begin{eqnarray}\label{eq:np}
n(F) &=&  \WH{\mc V}(x)^{-1} \, E\le[K_F(x)\ri],
\qquad
n^{(\gamma)}(F) =  \WH{\mc V}^{(\gamma)}({x})^{-1} E^{(\gamma)}\le[K_F({x})\ri].
\end{eqnarray}
In particular 
we have
$\WH{\mc V}(x)\cdot n(\rho(x) < \zeta)=1$
and
$\WH{\mc V}^{(\gamma)}(x)\cdot n^{(\gamma)}(\rho(x) < \zeta)=1$.

(ii) The following holds
$n^{(\gamma)}(F(\ve)\I_{\{\rho(x,\ve) < \zeta(\ve)\}}) 
=  n(\te{\gamma \ve(\rho(x,\ve))}F(\ve)\I_{\{\rho(x,\ve) < \zeta(\ve)\}})$.
Hence we have
\begin{equation}\label{eq:ng}
n^{(\gamma)}(\rho(x,\ve) < \zeta(\ve)) =  
n(\te{\gamma \ve(\rho(x,\ve))}\I_{\{\rho(x,\ve) < \zeta(\ve)\}}).
\end{equation}

(iii) For any $z\in(0,\infty)$ the following holds as $x\to\infty$:
\begin{equation}\label{eq:limitsnn}
n^{(\gamma)}(\rho(x,\ve) < \zeta(\ve))\simeq \WH\phi(\gamma)
\qquad\text{ and }\qquad
\te{\gamma x}n(\ve(\rho(z,\ve))>x, \rho(z,\ve)<\zeta(\ve)) = o(1).
\end{equation}
\end{Lemma}

\smallskip

\noindent {\it Proof of Lemma \ref{lem:npp}.}
(i) The proof of~\eqref{eq:np} is identical 
under both measures. Hence we give the argument 
only under 
$P$.
Note that for any left-end point 
$g$
of an excursion interval the following 
equality holds: 
$F(\e_g)\I_{\{g<\tau(x)\}} = F(\e_g)\I_{\{g\leq\tau(x)\}}$.
Since for every
$\ve\in\mc E$
the process
$t\to F(\ve)\I_{\{t\leq\tau(x)\}}$ 
is left-continuous and adapted,
an application of the compensation
formula of excursion theory for the Poisson point process
$\e$
defined in~\eqref{eq:e}
to $K_F(x)$ 
(see e.g.~\cite[Cor. IV.11]{Bertoin})
yields representation~\eqref{eq:np}.
The second statement  follows
by taking
$F=\I_{\{\rho(x)<\zeta\}}$
in~\eqref{eq:np}, since in that case $K_F(x) = \I_{\{\tau(x)<\infty\}}$.

(ii) Define
$G(\varepsilon)\doteq F(\ve)\I_{\{\rho(x,\ve) < \zeta(\ve)\}}$
and let
$K_G(x)$
as in~\eqref{eq:Def_K_F}.
The Esscher change of measure formula and 
the compensation formula 
in~\cite[Cor. IV.11]{Bertoin}
yield
\begin{eqnarray}
E^{(\gamma)}\le[K_G({x})\ri] 
\label{2}
&=& E\le[\int_0^\infty \te{\gamma X(t-)}\I_{\{t \leq \tau(x)\}}\td\WH L(t)\ri]
n\left(\te{\gamma\ve({\rho(x,\ve)})}F(\ve)\I_{\{\rho(x,\ve) < \zeta(\ve)\}}\right).
\end{eqnarray}
A change of variable
$t=\WH L^{-1}(u)$
under the expectation on the right-hand side of~\eqref{2},
Fubini's theorem
and
$P^{(\gamma)}$-a.s.
equality
$\{\WH L^{-1}(u-) \leq \tau(x)\}
=\{\WH L^{-1}(u) \leq \tau(x)\}$
yield
\begin{eqnarray*}
E\le[\int_0^\infty \te{\gamma X(t-)}\I_{\{t \leq \tau(x)\}}\td\WH L(t)\ri]
 = E^{(\gamma)}\le[\int_0^{\WH L(\infty)} \I_{\{\WH L^{-1}(u-) \leq \tau(x)\}} \td u\ri]
 =\WH{\mc V}^{(\gamma)}({x}).
\end{eqnarray*}
The final equality follows from
$\{\WH L^{-1}(u-) \leq \tau(x)\}
=
\{u \leq \WH L(\tau(x))\}$.
Equality in~\eqref{eq:np} under 
$P^{(\gamma)}$ applied to 
$K_G(x)$
and~\eqref{2}
now imply the formula in part~(ii) of the lemma.

(iii) By Lemma~\ref{lem:np} the random variable $\WH L({\tau(x)})$
is exponentially distributed under $P^{(\gamma)}$ with parameter
$n^{(\gamma)}(\rho(x)<\zeta)$. 
Hence
$n^{(\gamma)}(\rho(x)<\zeta) = -\log P^{(\gamma)}(\WH L({\tau(x)}) >1)$ 
and the dominated
convergence theorem implies
$\lim_{x\uparrow\infty} n^{(\gamma)}(\rho(x)<\zeta)=
- \log P^{(\gamma)}(\WH L(\infty) > 1) =
- \log P^{(\gamma)}(\WH L^{-1}(1) < \infty)$,
which is equal to
$\WH\phi^{(\gamma)}(0)= \WH\phi(\gamma)$
by the elementary equality
$\WH\phi^{(\gamma)}(u)= \WH\phi(\gamma +u)$,
$u\geq0$.
Chebyshev's inequality and part (ii) of the lemma imply
$\te{\gamma x}n(\ve(\rho(z,\ve))>x, \rho(z,\ve)<\zeta(\ve)) \leq 
n(\te{\gamma \ve(\rho(z,\ve))}\I_{\{\ve(\rho(z,\ve))>x,\rho(z,\ve)<\zeta(\ve)\}})
 = 
n^{(\gamma)}(\ve(\rho(z,\ve))>x,\rho(z,\ve)<\zeta(\ve)).$
The final expression tends to zero as $x\uparrow\infty$ by the dominated convergence theorem
and the lemma follows.
\exit

\medskip

We now apply Lemma~\ref{lem:npp}
to establish the asymptotic behaviour
of certain integrals against
the excursion measure
as
$x\to\infty$.
Lemma~\ref{lem:np}, 
in combination with Proposition~\ref{lem:exc} 
below, implies the identity in~\eqref{eq:zinf}
thus concluding the proof of 
Theorem~\ref{thm}. 

\begin{Prop}
\label{lem:exc}
Let $u\ge 0$.  Then, as $x\to\infty$, we have
\begin{equation}\label{eq:conle}
n(\te{-u(\ve({\rho(x)})-x)}|\rho(x)<\zeta) \longrightarrow C(u) \cdot C_\gamma^{-1}  =
\frac{\gamma}{\gamma + u}\cdot \frac{\phi(u)}{\phi(0)}.
\end{equation}
\end{Prop}


\noindent {\it Remark.}
Recall 
the result of Doney \& Maller~\cite[Thm.~1]{DoneyMaller} 
($C_\gamma$
is defined in~\eqref{eq:cra_const}):
\begin{equation}
\label{eq:DaneyMaller}
n(\rho(x) < \zeta) \simeq C_\gamma\, \WH\phi(\gamma)\, \te{-\gamma x}\qquad\text{as}\quad x\to\infty.
\end{equation}

\noindent {\it Proof of Proposition \ref{lem:exc}.}
Fix
$M>0$
and recall that, under the probability measure $n(\,\cdot\,|\rho(M)<\zeta)$,
the coordinate process has the same law as the first excursion
of $Y$ away from zero with height larger than $M$.
For any
$x>M$,
the following identity holds:
\begin{eqnarray}
\label{eq:Trivial_Useful}
n(\te{-u(\ve({\rho(x)})-x)}|\rho(x)<\zeta) & = &
n(\te{-u(\ve({\rho(x)})-x)}\I_{\{\rho(x)<\zeta\}} |\rho(M)<\zeta) \frac{n(\rho(M)<\zeta)}{n(\rho(x)<\zeta)}.
\end{eqnarray}

The definitions of the point process
$\epsilon$
in~\eqref{eq:e}
and of the compensator measure
$n$,
together with the strong Markov property
under the probability measure
$n(\,\cdot\,|\rho(M)<\zeta)$,
imply that
$\ve\circ\theta_{\rho(M)}$
has the same law as the process
$X$
with entrance law
$n(\ve({\rho(M,\ve)})\in\td z|\rho(M)<\zeta)$
and killed at the epoch of the first passage into
the interval
$(-\infty,0]$.
We therefore find
\begin{eqnarray}
\nonumber
\lefteqn{n(\te{-u(\ve({\rho(x,\ve)})-x)}\I_{\{\rho(x)<\zeta\}} |\rho(M)<\zeta)
= n\le(\te{-u(\ve({\rho(M,\ve)})-x)}\I_{\{\ve(\rho(M,\ve))> x\}}|\rho(M)<\zeta\ri)}\\
\label{eq:First_Markov}
& + & 
 \int_{[M,x]}
E_z\le[\te{- u K(x)}
\I_{\{T(x) < \WH T(0)\}}\ri]
n(\ve({\rho(M,\ve)})\in\td z|\rho(M)<\zeta),
\end{eqnarray}
where
$K(x)= X(T(x))-x$.
By the second equality in~\eqref{eq:limitsnn}
of Lemma~\ref{lem:npp}, 
we have 
as
$x\uparrow\infty$:
\begin{eqnarray*}
\te{\gamma x}n\le(\te{-u(\ve({\rho(M,\ve)})-x)}\I_{\{\ve(\rho(M,\ve))> x\}}|\rho(M)<\zeta\ri)
 \leq  \te{\gamma x}\frac{n\le(\ve(\rho(M,\ve))> x,\rho(M,\ve)<\zeta(\ve) \ri)}{n(\rho(M)<\zeta)}
 =  o(1).
\end{eqnarray*}
This estimate, 
spatial homogeneity of 
$X$
and equations~\eqref{eq:Trivial_Useful} and~\eqref{eq:First_Markov}
yield
as
$x\to\infty$:
\begin{eqnarray}
\nonumber
\lefteqn{n(\te{-u(\ve({\rho(x,\ve)})-x)}|\rho(x)<\zeta)}\\
&  = &  o(1) +
\int_{[M,x]}
E\le[\te{- u K(x-z)} \I_{\{T(x-z) < \WH T(z)\}}\ri]
\frac{n(\ve({\rho(M,\ve)})\in\td z, \rho(M)<\zeta)}{n(\rho(x)<\zeta)}.
\label{eq:Almost_Final}
\end{eqnarray}

Formula~\eqref{eq:Most_Important}
of Proposition~\ref{prop:aso}
implies the following equality:
\begin{equation}
E\le[\te{- u K(x-z)} \I_{\{T(x-z) < \WH T(z)\}}\ri] =
C(u) \te{-\gamma x} \left(1-G(z) +R(x-z)\right)\te{\gamma z}, 
\end{equation}
where
$G,R:\RR_+\to\RR$
are bounded functions such that 
$G(z) = E[\te{\gamma X(\WH T(z))}]$
and 
$\lim_{x'\to\infty}R(x')=0$.
Therefore the equality in~\eqref{eq:Almost_Final},
the asymptotic behaviour of
$n(\rho(x)<\zeta)$
given in~\eqref{eq:DaneyMaller}
and Lemma~\ref{lem:npp}~(ii)
imply the following identity
as
$x\to\infty$:
\begin{eqnarray}
\label{eq:Important_Expression}
\lefteqn{n(\te{-u(\ve({\rho(x,\ve)})-x)}|\rho(x)<\zeta)
=A_\gamma(u) n^{(\gamma)}(\ve({\rho(M,\ve)})\in[M,x],\rho(M,\ve)<\zeta(\ve)) +
o(1)} \\
& + &
A_\gamma(u)
 n^{(\gamma)}\left(\left[R(x-\ve({\rho(M,\ve)}))- G(\ve({\rho(M,\ve)}))\right]I_{\{\ve({\rho(M,\ve)})
 \in[M,x],\rho(M,\ve)<\zeta(\ve)\}}\right),
\nonumber
\end{eqnarray}
where
$A_\gamma(u)
\doteq
C(u)/(C_\gamma \WH \phi(\gamma))$.
By~\eqref{eq:Important_Expression}
the limit
$\lim_{x\to\infty}n(\te{-u(\ve({\rho(x,\ve)})-x)}|\rho(x)<\zeta)$
exists and the dominated convergence theorem yields 
\begin{eqnarray*}
\lim_{x\to\infty}n(\te{-u(\ve({\rho(x,\ve)})-x)}|\rho(x)<\zeta)
= A_\gamma(u) 
\left(n^{(\gamma)}(\rho(M)<\zeta) -
n^{(\gamma)}\left(G(\ve({\rho(M,\ve)}))I_{\{\rho(M,\ve)<\zeta(\ve)\}}\right)\right).
\end{eqnarray*}
Since this equality holds for any
$M>0$
and the left-hand side does not depend on
$M$,
if the right-hand side has a limit as
$M\to\infty$,
then the equality also holds in this limit.
Note that~\eqref{eq:limitsnn}
of Lemma~\ref{lem:npp}~(iii)
implies
$\lim_{M\to\infty}n^{(\gamma)}(\rho(M)<\zeta)= \WH \phi(\gamma)$.
Since 
$G(z) = E[\te{\gamma X(\WH T(z))}]$
it holds
$G(\ve({\rho(M,\ve)}))\leq\te{-\gamma M}$
and an application of the dominated convergence
theorem yields~\eqref{eq:conle}.
\exit


\appendix
\section{Additional proofs}

\subsection{Proof of Corollary~\ref{cor:4}.} 
\label{subsec:Cor}
(i) The duality lemma for L\'evy processes
implies that 
$X^*(t)=\sup_{0\leq s\leq t}X(s)$ 
and
$Y(t)$
have the same law for any fixed $t\geq0$.
Since, by As.~\ref{A2},
$E[X(1)]<0$
and the process  $\{X^*(t)\}_{t\ge 0}$
is non-decreasing, it
converges a.s. as $t\uparrow\infty$
to $X^*(\infty)\doteq\sup_{s\ge 0}X(s)$.  
Therefore
$Y(t)$ converges weakly to 
the law
$Y(\infty)$ 
of 
$X^*(\infty)$,
characterised by its Laplace transform
$E[\te{-u Y(\infty)}] = \phi(0)/\phi(u)$,
$u\in\RR_+$
(see~\cite[p.~163]{Bertoin}).
The joint Laplace transform 
of
$(Y(\infty),Z(\infty))$
now follows from the asymptotic independence in Theorem~\ref{thm}.

\noindent (ii) The Wiener-Hopf
factorisation of $X$ \cite[p. 166]{Bertoin}
implies the following identity 
for some $k\in(0,\infty)$: 
\begin{equation}\label{eq:WHF}
-\log E[\te{\theta X(1)}] = k \phi(-\theta) \WH\phi(\theta),
\qquad \theta\in\mathbb C,\ \Re(\theta) = 0.
\end{equation}
By analytic continuation and As.~\ref{A2}, 
identity~\eqref{eq:WHF} holds for all 
$\theta\in\mbb C$ 
with $\Re(\theta)\in[0,\gamma)$. 
Furthermore, continuity implies that~\eqref{eq:WHF} 
remains valid for 
$\theta=\gamma$.
As $\WH H$ is a non-zero subordinator (recall $E[X(1)]<0$),
we have
$\WH \phi(\gamma)>0$ 
and hence
$\phi(-\gamma)=0$.

By Thm.~\ref{thm}, eq.~\eqref{eq:zinf},
the Laplace transform of 
$x\mapsto P(Z(\infty) > x)$
is 
$(1-\frac{\gamma}{\phi(0)}\phi(v)/(v+\gamma))/v$.
The L\'evy-Khinchin formula for 
$\phi$
and integration by parts
imply
$\phi(v)=\phi(0)+v(m+\int_0^\infty\te{-vx}\ovl \nu_{H}(x)\,\td x)$
for any $v\geq-\gamma$.
Since 
$\phi(-\gamma)=0$,
we have
$ \int_0^\infty\te{\gamma y}\ovl\nu_H(y)\td y= \phi(0)/\gamma-m$.
A direct Laplace inversion, based on this representation of
$\phi$,
implies the first formula in part~(ii) of the corollary. 
The atom at zero 
is obtained by taking the limit
in~\eqref{eq:zinf}
of Theorem~\ref{thm}
as $v\to\infty$.

\noindent (iii) 
Since
$\phi$
is strictly concave
with 
$\phi(-\gamma)=0$,
the right-derivative
of $\phi$ 
at $-\gamma$
satisfies
$\phi'(-\gamma)>0$
and the constant 
$C_\gamma$
in~\eqref{eq:cra_const}
is well-defined.
By Lemma~\ref{lem:straddle}(i),
proved in Appendix~\ref{subsec:Local_time_expectation},
we have 
$\ell\in(0,\infty)$.
It follows from~\cite[Thm.~1]{DoneyMaller}
that if $t$ and $x$ tend to infinity and
$t\te{-\gamma x}\rightarrow \lambda$, for some $\lambda>0$,
then $M({t,x})\stackrel{\mc D}{\longrightarrow} M(\infty)$, where
the limit follows a Gumbel distribution,
\begin{equation}
\label{eq:limyt} P(M(\infty) < z) = \exp\le(-\ell\, C_\gamma\, \WH\phi(\gamma)\,
\lambda\,\te{-\gamma z}\ri),\qquad \text{for all}\quad z\in\RR.
\end{equation}
For completeness, we give below a short proof of~\eqref{eq:limyt}
based on~\cite[Thm.~1]{DoneyMaller}. 
The joint Fourier-Laplace transform in Corollary~\ref{cor:4}
now follows from Theorem~\ref{thm} and a direct calculation using~\eqref{eq:limyt}.

To establish~\eqref{eq:limyt}
we show that, as 
$t\wedge x\to\infty$ and $t\te{-\gamma x}\to\lambda>0$, the
following holds
\begin{eqnarray}
P(Y^*(t)-x < z) &=& 
\exp(-t\,\ell\, \lambda(x+z)) + o(1)\quad\text{for any $z\in\RR$},
\label{eq:Ystar}
\end{eqnarray}
where $\lambda(x)= n(\rho(x)<\zeta)$. Since~\eqref{eq:DaneyMaller}
implies $t\lambda(x+z)\to C_\gamma\WH \phi(\gamma)\lambda\te{-\gamma z}$ as
$t\wedge x\to\infty$ and $t\te{-\gamma x}\to\lambda$,
the limit in~\eqref{eq:limyt} follows
from~\eqref{eq:Ystar}.

To complete the proof we now verify the claim in~\eqref{eq:Ystar}.
Note that $\tau(x+z)\to\infty$ $P$-a.s. as $x\to\infty$ and, as
shown in the proof of Lemma~\ref{lem:straddle}, the law of large
numbers implies that $\WH L(t)/t\to\ell$ $P$-a.s. as $t\to\infty$,
where $\ell=1/E\left[\WH L^{-1}(1)\right]$ (recall from
Lemma~\ref{lem:straddle}~(i) that $0<\ell<\infty$). Therefore 
$\WH L(\tau(x+z))/\tau(x+z)$ 
tends to $\ell$ $P$-a.s as $x\to\infty$. In
particular, for any $\delta>0$, we have
$$P(\WH L(\tau(x+z))/\tau(x+z)\in(\ell-\d,\ell+\d)) = 1 + o(1)\qquad\text{as}\quad x\to\infty.$$
Hence
as $t\wedge x\to\infty$ the following holds
\begin{eqnarray*}
P(Y^*(t) < x+z) & = & P(\tau(x+z) > t, \WH L(\tau(x+z))/\tau(x+z)\ge \ell-\d)  + o(1)\\
&\leq& P(\WH L({\tau(x+z)}) > t(\ell-\d)) + o(1).
\end{eqnarray*}
Similarly, it follows that as $t\wedge x\to\infty$ we have
\begin{eqnarray*}
P(Y^*(t) < x+z) & \ge & P(\WH L({\tau(x+z)})> \WH L(t), \WH L(t)\leq t(\ell+\delta))\\
&\ge& P(\WH L({\tau(x+z)})>t(\ell+\d), \WH L(t)\leq t(\ell+\delta)) = P(\WH L({\tau(x+z)})>t(\ell+\d)) + o(1).
\end{eqnarray*}
By Lemma~\ref{lem:np} the random variable $\WH L(\tau(x+z))$ is
exponentially distributed with parameter $\lambda(x)$ and hence we
find
\begin{equation*}
\exp(-(\ell+\d)t\, \lambda(x+z)) + o(1) \leq P(Y^*(t) < x+z)  \leq
\exp(-(\ell-\d)t\,\lambda(x+z)) + o(1).
\end{equation*}
Since this result holds for any $\delta>0$, the equality in
\eqref{eq:Ystar} follows.\exit

\subsection{Proof of Proposition~\ref{prop:aso}} (i) 
\label{subsec:Two_sided_Overshoot_LEvy}
Recall that 
$T(x)$ 
and
$\WH T(x)$ 
are defined in~\eqref{eq:TT}
and that, 
under As.~\ref{A2},
\cite{BertoinDoney} shows that
Cram{\'e}r's estimate remains valid for the  L\'{e}vy process $X$
($C_\g$
defined in~\eqref{eq:cra_const}):
\begin{equation}\label{eq:cra}
P(T(y)<\infty)
\simeq C_\g\te{-\gamma y}
\quad \text{as $y\to\infty$}.
\end{equation}
By the
strong Markov property and spatial homogeneity of $X$ it follows
that
\begin{equation}
\label{eq:Markov_Prop_Simple_Case} P(T(x) < \WH T(z))
= P(T(x)<\infty) - \int_{(-\infty,-z]}P_y(T(x)<\infty) P(X({\WH
T(z)})\in\td y, \WH T(z)<T(x)).
\end{equation}
The translation invariance of $X$ and Cram\'er's
estimate~\eqref{eq:cra} imply the following equality
\begin{eqnarray}
P_y(T(x)<\infty) 
& = & C_\gamma \te{-\gamma x} \te{\gamma y}
\left(1+r(x-y)\right)\quad\text{for all}\quad x>y,
\label{eq:Useful_Two_sided_Asym}
\end{eqnarray}
where $\lim_{x'\to\infty} r(x') = 0$.
Equality~\eqref{eq:Useful_Two_sided_Asym} applied to the identity
in~\eqref{eq:Markov_Prop_Simple_Case} yields
\begin{eqnarray}
\label{eq:key_identity} C_\gamma^{-1} \te{\gamma x} P(T(x) < \WH
T(z)) & = &
1 - E\le[\te{\gamma X({\WH T(z)})}\I_{\{\WH T(z)<T(x)\}}\ri]\\
&+& r(x) - E\le[\te{\gamma X(\WH T(z)) }r(x-X(\WH T(z)))\I_{\{\WH
T(z)<T(x)\}}\ri]. \nonumber
\end{eqnarray}
Since $X({\WH T(z)})\leq -z<0$ on the event $\{\WH T(z)<\infty\}$,
which satisfies $P(\WH T(z)<\infty)=1$ by As.~\ref{A2},
the dominated convergence theorem implies
$$
E\le[\te{\gamma X({\WH T(z)})}\ri] = 
E\le[\te{\gamma X({\WH T(z)})}\I_{\{\WH T(z)<T(x)\}}\ri]  +
o(1)\qquad\text{as $x\to\infty$.}
$$
An application of the dominated convergence theorem to the second
expectation on the right-hand side of
equality~\eqref{eq:key_identity}, together with the fact that $r$
vanishes in the limit as $x\to\infty$, proves the first statement
in the proposition. 

(ii) Recall that the Laplace exponent $\phi$ of the increasing ladder-height
process $H$ is a strictly concave function that satisfies
$\phi(-\gamma)=0$ so that the right-derivative $\phi'(-\gamma)$ is
strictly positive. Under the measure $P^{(\gamma)}$ the identity
$\phi^{(\gamma)}(\gamma + u)=\phi(u)$ holds for any $u\in\mbb R_+$
and hence, since $\phi'(-\gamma) = E^{(\gamma)}[X_1]>0$, $X$
drifts to $+\infty$ as $t\to\infty$, i.e.
$P^{(\gamma)}(T(x)<\infty)=1$ 
for any
$x>0$.
Therefore, under As.~\ref{A2}, 
under 
$P^{(\gamma)}$ 
the ladder-height process 
$H$
is a non-lattice subordinator with
$E^{(\gamma)}\left[H(1)\right]\in(0,\infty)$.
Since the overshoot 
$K(x)$
is equal to that of
$H$
over $x$,
\cite[Thm.~1]{BertoinvanHarnSteutel}
implies that the weak limit 
$K(x) \stackrel{\mc D}{\longrightarrow} K(\infty)$,
as $x\to\infty$,
exists.
Since
$x\mapsto \te{-ux}$
is uniformly continuous on $\RR_+$, 
\cite[p.~16,~Thm.~2.1]{Billingsley}
implies
$\lim_{x\uparrow\infty}E^{(\gamma)}[\te{-uK(x)}] =
E^{(\gamma)}[\te{-uK(\infty)}]$
for any fixed $u\geq0$.
A version of the Wiener-Hopf factorisation of $X$ (see e.g. \cite[p.183]{Bertoin})
under the measure 
$P^{(\gamma)}$
yields
\begin{equation}\label{eq:LT}
\int_0^\infty q\te{-qx} E^{(\gamma)}\left[\te{-u K(x)} \right]\>\td x = \frac{q}{\phi(q-\gamma)}\cdot\frac{\phi(q-\gamma)-\phi(u-\gamma)}{q-u}\qquad
\text{for any $q,u>0$.}
\end{equation}
Since  the function 
$x\mapsto E^{(\gamma)}\left[\te{-u K(x)} \right]$
is bounded, the dominated convergence theorem implies that
in the limit as 
$q\downarrow0$
we get
$E^{(\gamma)}[\te{-uK(\infty)}]= \phi(u-\gamma)/(u\phi^\prime(-\gamma))$. 
The Esscher change of measure formula 
implies the following for 
any $u\ge 0$
($C(u)$ is defined in~\eqref{eq:Most_Important}): 
\begin{equation}
\label{eq:chme} E[\te{-u K(x)}\I_{\{T(x)<\infty\}}] = \te{-\gamma x} \cdot E^{(\gamma)}[\te{-(\gamma+u)K(x)}]
\simeq C(u) \te{-\gamma x} \quad \text{as $x\to\infty$}. 
\end{equation}
Furthermore, since the expectation in~\eqref{eq:chme} is bounded
as $x\to\infty$, there exists a bounded function $R:\RR_+\to\RR$,
such that
$E[\te{-u K(x)}\I_{\{T(x)<\infty\}}] = C(u) \te{-\gamma x}(1+R(x))$
for $x>0$, and
$\lim_{x\to\infty}R(x)=0$.
The strong Markov property at $\WH T(z)$
and an argument analogous to the one used in the proof of
Proposition~\ref{prop:aso}(i)
(cf.~\eqref{eq:key_identity}) yields
\begin{eqnarray*}
\lefteqn{C(u)^{-1}\te{\gamma x}E[\te{-u K(x)}\I_{\{T(x)<\WH T(z)\}}]} \\
& = &   1  - E[\te{\gamma X(\WH T(z))} \I_{\{\WH T(z) < T(x)\}}] +
R(x) - E[\te{\gamma X(\WH T(z))}R(x-X(\WH T(z)))\I_{\{\WH T(z) <
T(x)\}}],
\end{eqnarray*}
which implies equivalence~\eqref{eq:Most_Important}.\exit

\subsection{Proof of Lemma~\eqref{lem:straddle}(i)} 
\label{subsec:Local_time_expectation}
By analytical continuation and As.~\ref{A2} it follows
that identity~\eqref{eq:WHF} remains valid for all $\theta\in\mbb
C$ with $\Re(\theta)\in[0,\gamma)$. Therefore on the event
${\{H(1)<\infty\}}$ the random variable $H(1)$ admits finite
exponential moments and in particular
$E\left[H(1)\I_{\{H(1)<\infty\}}\right]<\infty$. Since
$E[X(1)]\in(-\infty,0)$, the ladder-height process of the dual
process $\WH X=-X$ satisfies $P(\WH H(1)<\infty)=1$. 
Furthermore, we
have $P(H(1)<\infty)<1$. Definition~\eqref{eq:phi} of $\phi$, its
analogue for $\WH \phi$, the Wiener-Hopf factorisation
in~\eqref{eq:WHF} and the dominated convergence theorem imply
that the following identity holds for all $\theta\in(0,\gamma)$:
\begin{eqnarray*}
-\frac{E[X(1)\te{\theta X(1)}]}{kE[\te{\theta X(1)}]} & = &
\frac{E[H(1)\te{\theta H(1)}\I_{\{H(1)<\infty\}} ]}{E[\te{\theta
H(1)}\I_{\{H(1)<\infty\}} ]}
\log E\left[\te{-\theta \WH H(1)}\right]\\
& - & \frac{E[\WH H(1)\te{-\theta \WH H(1)}]}{E[\te{-\theta \WH
H(1)}]} \log E\left[\te{\theta H(1)}\I_{\{H(1)<\infty\}} \right].
\end{eqnarray*}
As.~\ref{A2}
implies that in the limit as $\theta\to0$ this equality yields
$E[\WH H(1)]\in(0,\infty)$.

The inequality $ \left|\WH X\left(t\wedge \WH
L^{-1}(1)\right)\right|\leq \WH H(1) + X^*(\infty)$ holds for all
$t\in\RR_+$. 
Cram\'er's estimate~\eqref{eq:cra} implies that $X^*(\infty)$ is integrable.
Since $\left\{\WH X(t)-t E[\WH X(1)]\right\}_{t\geq0}$ 
is a martingale we have
$$E\left[\WH X\left(t\wedge \WH L^{-1}(1)\right)\right]  =  E\left[\WH X(1)\right]
E\left[t\wedge \WH L^{-1}(1)\right]\quad\text{for all}\quad
t\in\RR_+.$$ The dominated and monotone convergence theorems
applied to each side of this equality respectively imply Wald's
identity for the $\{\mathcal F_t\}$-stopping time 
$\WH L^{-1}(1)$:
$E\left[\WH H(1)\right]  =  - E\left[X(1)\right]E\left[\WH L^{-1}(1)\right].$
In particular we obtain
$\ell^{-1} = E\left[\WH L^{-1}(1)\right]\in(0,\infty),$
proving Lemma~\eqref{lem:straddle}(i).\exit

\bibliographystyle{plain}
\bibliography{cite}

\end{document}